\documentclass[amsmath,english,a4paper,graphicx,12pt]{article}
\usepackage{doc}
\usepackage{graphicx}
\usepackage{amsmath}
\usepackage{amsfonts}
\usepackage{amssymb}
\usepackage{babel}
\usepackage{latexsym}
\usepackage{mathrsfs}

\title{\bf Development of a complex function theory upon a  new concept of polar-analytic functions;\\  Extended version}
\author{Carlo Bardaro, \thanks{
Department of Mathematics and Computer Sciences, University of Perugia,
via Vanvitelli 1, I-06123 Perugia, Italy, e-mail: 
carlo.bardaro@unipg.it} \and
Paul L. Butzer, \thanks{Lehrstuhl A fuer Mathematik, RWTH Aachen, Templergraben 55, Aachen, D-52056, Germany, e-mail: butzer@rwth-aachen.de}\and
Ilaria Mantellini, \thanks{Department of Mathematics and Computer Sciences, University of Perugia,
via Vanvitelli 1, I-06123 Perugia, Italy, e-mail: 
mantell@dmi.unipg.it}\and Gerhard Schmeisser\thanks{Department Mathematik, FAU Erlangen-N\"{u}rnberg, Cauerstr. 11, 91058 Erlangen, Germany, 
email: schmeisser@mi.uni-erlangen.de}
}
\begin{document}
\maketitle
\noindent

\begin{abstract}
The present article is an extended version of  \cite{BBMS4} containing new results and an updated list of references. We review the notion of polar analyticity introduced in \cite{BBMS1} and succesfully applied in Mellin analysis and quadrature formulae for functions defined on the positive real axis. This appears as a simple way to describe functions which are analytic on a part of the Riemann surface of the logarithm. We also describe some geometric properties of polar-analytic functions related to conformality.
In this paper we launch a proposal to develop a complete complex function theory, independent of classical function theory, which is built upon the new concept of polar analyticity.

\end{abstract}
%-----------------------------------------------------------------------------------------
%------------------------INTRODUZIONE----------------------------------------------------

%------------------------Preliminari----------------------------------------------------

%---------------------------------------------- --------------------------------

%-------------------------------------------section----------------------------------
\section{The background}

Let us first recall that in the classical complex analysis the Cauchy-Riemann equations  in polar coordinates are obtained by setting $z= re^{i \theta},$ i.e., given an analytic function $f= u + iv$ on a domain in the complex plane, one considers the function $g(r,\theta) := f(re^{i\theta})$ and uses the chain rule for deriving partial derivatives of $u,v$  with respect to the variables $r,\theta,$ (see e.g. \cite[Sec. 23, p. 68]{BC}; \cite[Sec. 4.3, pp. 78--82]{HI}). Due to the $2\pi$-periodicity of the exponential function $e^{i\theta},$ this implies a periodicity with respect to $\theta$ of the function $g.$ As we will see, in our definition of polar analyticity we again derive the same Cauchy-Riemann equations in polar form, but in general this periodicity does not appear.

In order to define the notion of polar-analytic function on a domain in the polar plane, we begin with some preliminary facts, which are in some sense connected with the periodicity with respect to $\theta.$
Corresponding to each nonzero complex number $z = r e^{i\theta},$ the multivalued function (complex logarithm)
$$\log z = \log r + i \theta$$
can be described as a single-valued function by replacing the complex plane without the origin, by the so-called Riemann surface of the logarithm (see e.g. \cite[Sec. 4.3, p. 97]{AH}). This is defined in an abstract way as a connected surface of infinitely many sheets, and on each sheet $\theta$ ranges in an interval of width $2 \pi.$ 

Let us denote by $S_{{\rm log}}$ the Riemann surface of the complex logarithm. 
A simple and convenient model of $S_{{\rm log}}$ is the helicoidal surface in $\mathbb{R}^3$ defined by
$$ \mathbb{E}:= \{ (x_1,x_2,x_3) \in \mathbb{R}^3 : x_1= r \cos \theta,~ x_2= r\sin \theta,~x_3=\theta,~ r>0, ~\theta\in \mathbb{R}\}.$$
The subset obtained by setting $\theta =0$ on the right-hand side can be interpreted as $\mathbb{R}^+.$
Just as $\mathbb{C}$ is an extension of  $\mathbb{R},$ we shall see that  $\mathbb{E}$ takes a corresponding role for the positive real axis 
$\mathbb{R}^+.$ This is very useful in Mellin analysis, when one wishes to extend a Mellin bandlimited function defined on $\mathbb{R}^+,$ to the 
Riemann surface  $S_{{\rm log}}.$
\begin{figure}[htbp] 
\begin{center}
\includegraphics[scale=1.00]{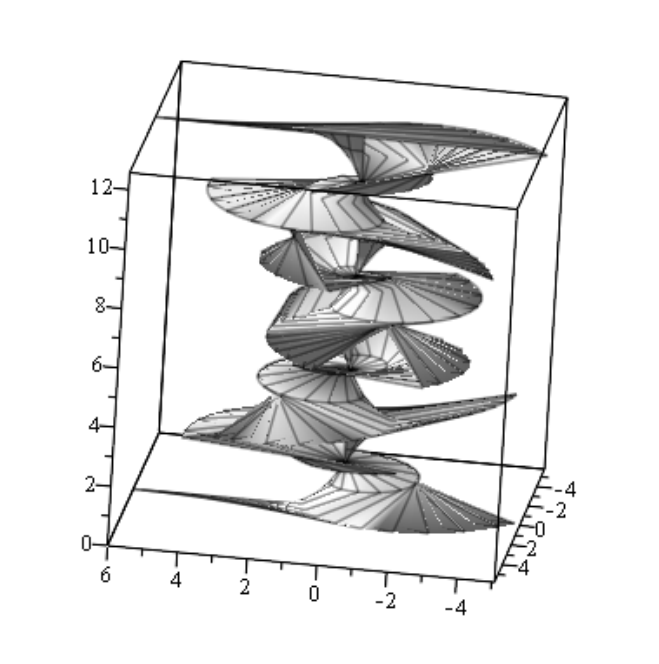}
\hskip0.2cm
\end{center}
\caption{\small The helicoidal surface as a model of the Riemann surface of the logarithm.}\label{helix}
\end{figure}
For $\alpha,~\beta\in \mathbb{R}$ with $\alpha<\beta,$ we consider the surface
$$\mathbb{E}_{\alpha,\beta}:= \left\{(x_1,x_2,x_3) \in \mathbb{R}^3 : x_1= r \cos \theta,~ x_2= r\sin \theta,~x_3=\theta,~ r>0, ~\theta\in {]}\alpha,\beta{[}\right\}$$
and call it a segment of $\mathbb{E}.$ The projection of $\mathbb{E}_{\alpha,\beta}$ into the $(x_1,x_2)$-plane may be interpreted as a sector in $\mathbb{C},$ given by
$$ S_{\alpha,\beta}:= \left\{ x_1 +ix_2 \in \mathbb{C} : x_1= r\cos\theta,~ x_2=r\sin \theta,~ r>0,~ \theta\in {]}\alpha,\beta{[} \right\}.$$
When $\beta -\alpha \in~{]}0, 2\pi],$ then this projection is a bijection. Indeed, given $z\in  S_{\alpha,\beta},$ there exists a unique $\theta\in~ {]}\alpha,\beta{[},$ denoted by $\theta:= \mbox{arg}_\alpha z,$ such that $z= |z| e^{i {\rm arg}_\alpha z}$ and then  $(\Re z, \Im z  ,\mbox{arg}_\alpha z)$ is the pre-image of $z$ on $\mathbb{E}_{\alpha,\beta}.$
\vskip0,3cm
\section{Polar-analytic functions}

Now, we can introduce analytic functions on $ \mathbb{E}$ as follows (see \cite{BBMS1}):
\vskip0,3cm

{\bf Definition 1}. 
A function $f: \mathbb{E} \rightarrow \mathbb{C}$ is said to be {\rm analytic} if for every segment 
$\mathbb{E}_{\alpha,\beta}$ with $\beta-\alpha\in~]0, 2\pi]$ the function $z\longmapsto f(\Re z, \Im z, {\rm arg}_\alpha z)$ is analytic on $S_{\alpha,\beta}.$
\vskip0,3cm
As an example, the function $L$ defined by
$$ L(r\cos \theta, r\sin \theta, \theta):= \log r+i \theta \quad \quad (r>0,~ \theta \in \mathbb{R}),$$
is analytic on $\mathbb{E}$ and coincides on $\mathbb{R}^+$ with the logarithm of real analysis.

Given $\alpha,$ the largest admissible $\beta$ in the previous definition is $\alpha + 2\pi.$
For $k\in \mathbb{Z},$ all the sectors $S_{\alpha+2k\pi,\alpha +2(k+1)\pi}$ coincide with the complex plane slit along the ray $z=r e^{i \alpha} (r>0)$ and the analytic functions induced by $f$ on these sectors are the analytic branches of $f.$

The above approach to analytic functions on  $S_{{\rm log}}$   is more understandable, due to the simple geometric representation of the Riemann surface. However, in this approach
we have functions of three variables $x_1, x_2, x_3$ but practically we use them only as functions of the 
two variables $r$ and $\theta.$ Since there exists a bijection between the helicoidal surface and the right half-plane understood as the set of 
all points $(r, \theta)$ with $r > 0$ and $\theta \in \mathbb{R},$ one may think of considering functions defined on the right half-plane. However, 
these functions will no longer be analytic in the classical sense. They are differentiable and satisfy the Cauchy--Riemann equations transformed 
into polar coordinates. This new approach amounts to taking an analytic function, writing its variable in polar coordinates $z = r e^{i \theta}$ and 
treating $(r, \theta)$ as if they were Cartesian coordinates (see \cite{BBMS1}, \cite{BBMS2}).

Let $\mathbb{H}:= \{(r,\theta) \in \mathbb{R}^+ \times \mathbb{R}\}$ be the right half-plane and let ${\cal D}$ be a domain in $\mathbb{H}.$
\vskip0,3cm
{\bf Definition 2}.
We say that $f:{\cal D}\rightarrow \mathbb{C}$ is {\rm polar-analytic} on ${\cal D}$ if for any $(r_0, \theta_0) \in {\cal D}$ the limit
$$\lim_{(r,\theta) \rightarrow (r_0, \theta_0)}\frac{f(r, \theta) - f(r_0, \theta_0)}{re^{i\theta} - r_0e^{i\theta_0}} =: (D_{{\rm pol}}f)(r_0, \theta_0)$$
exists and is the same howsoever $(r, \theta)$ approaches $(r_0, \theta_0)$ within ${\cal D}.$
\vskip0,3cm
\newtheorem{Remark}{Remark}
\begin{Remark}
{\rm  Let $f\,:\, (r,\theta) \mapsto u(r,\theta)+ iv(r,\theta)$ be
polar-analytic on  ${\cal D}.$ Then, identifying $\mathbb{C}$ with $\mathbb{R}^2$, we may interpret $f$
as a mapping from a subset of the half-plane $\mathbb{H}$ into $\mathbb{R}^2$. We note
that this mapping is differentiable in the classical sense. Indeed, polar 
analyticity implies that
\begin{equation}\label{1}
f(r,\theta)-f(r_0,\theta_0) = 
\left(D_{\text{pol}} f\right)(r_0,\theta_0)\left(re^{i\theta}-
r_0e^{i\theta_0}\right) + o\left(\left|re^{i\theta}- r_0e^{i\theta_0}\right|\right)
\end{equation}
as $re^{i\theta}\to r_0e^{i\theta_0}.$ But $re^{i\theta}$ is
differentiable as a mapping from $\mathbb{H}$ to $\mathbb{R}^2$. We have
\begin{equation}\label{2}
re^{i\theta}- r_0 e^{i\theta_0}= e^{i\theta_0}(r-r_0) + ir_0e^{i\theta_0}
(\theta-\theta_0)
+o\left(\left\|{r\choose\theta}-{r_0\choose\theta_0}\right\|_2\right)
\end{equation}
as $(r,\theta)\to(r_0,\theta_0)$, where $\| \cdot \|_2$ is the euclidean
norm for column vectors of $\mathbb{R}^2$.
Now, combining (\ref{1}) and (\ref{2}), and noting that
$$ \left(D_{\text{pol}}\right)(r,\theta) \,=\,
e^{-i\theta}\,\frac{\partial}{\partial r} f(r,\theta)\,=\,
e^{-i\theta}\,\left(\frac{\partial}{\partial r} u(r,\theta)
	+i\,\frac{\partial}{\partial r} v(r,\theta)\right),$$
we arrive at
\begin{align}
f(r,\theta)-f(r_0,\theta_0)\,&=\, 
{u(r,\theta)-u(r_0,\theta_0)\choose v(r,\theta)-v(r_0,\theta_0)} 
	\nonumber\\[2ex]
& =\, J(r_0,\theta_0) {r-r_0\choose \theta-\theta_0} +
o\left(\left\|{r\choose\theta}-{r_0\choose\theta_0}\right\|_2\right) \label{3}
\end{align}
as $(r,\theta)\to(r_0,\theta_0)$, where
$$ J(r_0,\theta_0)\,:=\, 
	\begin{pmatrix}
{\displaystyle\frac{\partial u}{\partial r}} & -r_0\,{\displaystyle\frac{\partial v}{\partial r}}
	\\[2ex]
{\displaystyle\frac{\partial v}{\partial r}} & \hphantom{-}r_0\,{\displaystyle\frac{\partial u}{\partial r}}
	\end{pmatrix}
$$
with the partial derivatives evaluated at $(r_0,\theta_0)$  takes the role
of the Jacobi matrix in case of polar-analytic functions. This shows that
$f$ interpreted as a mapping into $\mathbb{R}^2$ is differentiable in the
classical sense.}
\end{Remark}
\vskip0,3cm
\begin{Remark}
{\rm It can be verified that $f = u + iv$ with $u,v: {\cal D}\rightarrow \mathbb{R}$ is polar-analytic on ${\cal D}$ if and only if $u$ and $v$ have continuous partial derivatives on ${\cal D}$ that satisfy the differential equations}
\begin{eqnarray}\label{CRE}
\frac{\partial u}{\partial \theta} = - r \frac{\partial v}{\partial r}\,,\quad 
\frac{\partial v}{\partial \theta} = r \frac{\partial u}{\partial r}\,.
\end{eqnarray}
\end{Remark}

Note that these equations coincide with the Cauchy-Riemann equations of an analytic function $g$ defined by $g(z) := u(r,\theta) + i v(r, \theta)$ for $z= r e^{i\theta}.$ For the derivative $D_{{\rm pol}}$ we easily find that 
\begin{eqnarray}\label{dipol}
(D_{{\rm pol}}f)(r, \theta) = e^{-i\theta}\bigg[\frac{\partial}{\partial r}u(r, \theta) + i \frac{\partial}{\partial r}v (r, \theta) \bigg] = 
\frac{e^{-i\theta}}{r}\bigg[\frac{\partial}{\partial \theta}v (r, \theta) - i \frac{\partial}{\partial \theta}u (r, \theta) \bigg].
\end{eqnarray}
Since $f=u+iv,$ equations (\ref{CRE}) can be written in a more compact way as
\begin{eqnarray}\label{CRE2}
\frac{\partial f}{\partial \theta}  = ir \frac{\partial f}{\partial r}
\end{eqnarray}
and then formula (\ref{dipol}) takes the form
$$(D_{{\rm pol}}f)(r, \theta) = e^{-i\theta}\frac{\partial}{\partial r}f(r,\theta) = \frac{e^{-i\theta}}{ir} \frac{\partial}{\partial \theta}f(r,\theta).$$

Also note that $D_{{\rm pol}}$ is ordinary differentiation on $\mathbb{R}^+.$ More precisely, if $\varphi (\cdot) := f(\cdot, 0),$ then $(D_{{\rm pol}}f)(r,0) =
\varphi'(r).$

When $g$ is an entire function, then $f: (r, \theta) \mapsto g(re^{i\theta})$ defines a function $f$ on $\mathbb{H}$ that is polar-analytic 
and $2\pi$-periodic with respect to $\theta.$ Moreover by (\ref{dipol}) one has $(D_{{\rm pol}}f)(r, \theta) = g'(z)$ with  $z=re^{i\theta}.$ 

Conversely, if $f$ is polar-analytic on $\mathbb{H}$ and is $2\pi$-periodic, we cannot deduce in general the existence of an entire function $h$ such that $f(r,\theta) = h(re^{i\theta}).$ A simple example is the function $f(r,\theta):= e^{-i\theta}/r,$ for which we take $h(z) = 1/z$ that is analytic on $\mathbb{C}\setminus \{0\}.$
However, if $f$ is a polar-analytic function on $\mathbb{H},$ then $g: z=x+iy\longmapsto f(e^x, y)$ is an entire function (see Proposition \ref{taylor} below).

As examples, let us consider the function $g(z)= z^a,$ $a>0.$ Let us put $f(r,\theta):= g(re^{i\theta}) = r^ae^{ia\theta}.$ Then one has
$$f(r,\theta) = u(r,\theta) + i v(r,\theta):= r^a \cos (a\theta) + ir^a\sin(a\theta)$$
and so
$$(D_{{\rm pol}}f)(r, \theta) = e^{-i\theta}\bigg[\frac{\partial u}{\partial r}(r,\theta) + i\frac{\partial v}{\partial r}(r,\theta)\bigg] = a(re^{i\theta})^{a-1} = g'(z).$$
Analogously, putting $g(z) = \sin z$ and $f(r,\theta):= \sin(re^{i\theta}),$ we have $(D_{{\rm pol}}f)(r, \theta) = \cos(re^{i\theta}) = g'(z).$ 
Related examples may be found e.g. in \cite[Sec. 23]{BC} and \cite[Sec. 4.3]{HI}.

The main novelty of the definition of polar-analytic function is that using this approach we avoid periodicity with respect to the argument $\theta,$ and in this way we can avoid the use of Riemann surfaces. 

A simple example of polar-analytic function that is not $2\pi$-periodic is the function $L(r, \theta):= \log r + i\theta,$ which is easily seen 
to satisfy the differential equations (\ref{CRE}). In this approach, we consider the logarithm as a single-valued function on $\mathbb{H},$ without the use of the Riemann surface $S_{{\rm log}}.$ Using (\ref{dipol}) we find ($z=re^{i\theta}$)
$$(D_{{\rm pol}} L)(r, \theta) = e^{-i\theta}\frac{1}{r} = \frac{1}{re^{i\theta}} = \frac{1}{z}.$$
\vskip0,4cm
In order to state a connection with analytic functions on $S_{{\rm log}},$ 
for $\alpha, \beta \in \mathbb{R}$ with $\alpha < \beta,$ we consider the set
$$\mathbb{H}_{\alpha, \beta}:= \{(r, \theta) \in \mathbb{R}^+ \times \mathbb{R}: \theta \in ]\alpha, \beta[\}$$
and call it a {\it strip} of $\mathbb{H}.$

If $f : \mathbb{H} \rightarrow \mathbb{C}$ is polar-analytic but not $2\pi$-periodic with respect to $\theta,$ then we can associate with $f$ 
a function $g$ that is analytic on the Riemann surface $S_{{\rm log}}$ of the logarithm. The restriction of $f$ to a strip 
$\mathbb{H}_{\alpha + 2k\pi, \alpha + 2(k+1)\pi},$ where $k \in \mathbb{Z},$ defines an analytic function $g_k$ in the slit complex plane 
$\mathbb{C}\setminus \{re^{i\alpha}: r>0\}$ by setting $g_k(re^{i\theta}):= f(r,\theta).$ The functions $g_k$ for $k \in \mathbb{Z}$ are 
the analytic branches of $g.$
\vskip0,4cm
It seems to us that this modified notion of analyticity, arising by treating polar coordinates as 
Cartesian coordinates, has not yet been presented.   Our definition leads naturally to  the classical Cauchy-Riemann
equations when written in their  polar form, often treated
in the literature (see \cite{BC}, \cite{HI}). Although other mathematicians may have come across this concept too (see e.g. \cite{SI}), it seems that it has not been used for practical purposes so far. 

In Mellin analysis it turns out to be very helpful for an efficient approach, independent of Fourier analysis. In particular it leads to a precise and 
simple analysis for functions defined on the Riemann surface of the complex logarithm, via the helicoidal surface. 
 For details see \cite{BBMS1}, \cite{BBMS2}.

Different studies of analytic functions over the Riemann surface of the logarithm, identified with $\mathbb{H},$ were presented in papers dealing with Dirichlet's problem in certain domains, the Riemann mapping theorem, $o$-minimality and  the extension of analytic function on a domain with a boundary with prescribed properties (see \cite{KA1}, \cite{KA2}, \cite{KRS}). However, the notion of derivative introduced by these authors is different from ours (see \cite[Definition 2.1]{KRS}) and it is related to the Mellin derivative of a function (see \cite[Chapter 2, page 61]{MA}, \cite{BJ} and the subsequent Section 4).

The aim of this paper is to propose a possible development of a complete complex function theory, independent of classical function theory, which is built upon the new concept of polar analyticity.

 In order to state a version of Taylor's expansion theorem for polar-analytic functions, which establishes a fundamental connection with the analytic functions in the classical sense, we premise some remarks.

If $z_0 = x_0 + i y_0 \in \mathbb{C},$ for $\rho >0,$ let $D(z_0, \rho):= \{z \in \mathbb{C}: |z-z_0| < \rho\}$ be a disk in $\mathbb{C}.$ The above disk is transformed into the region with one axis of symmetry in $\mathbb{H}$ given by
$$E((r_0,\theta_0), \rho):= \left\{(r,\theta) \in \mathbb{H}: \bigg(\log	\frac{r}{r_0}\bigg)^2 + (\theta -\theta_0)^2 < \rho^2\right\},$$
which we call a "polar-disk" with "center" $(r_0,\theta_0)$ and "radius" $\rho.$. Here, $r_0= e^{x_0}$ and $\theta_0 = y_0.$ The boundary of this region is given by the graph of the functions
$$r = r_0\exp\big(\pm \sqrt{\rho^2 - (\theta-\theta_0)^2}\big), \quad \theta \in [-\rho + \theta_0, \rho + \theta_0].$$
Given  a domain ${\cal D} \subset \mathbb{H},$ let us define
$$A:= \{z=x+iy \in \mathbb{C}: (e^x,y) \in {\cal D}\}.$$
Then, for $z_0 = x_0+iy_0 \in A$ and a disk $D(z_0,\rho) \subset A,$ in order that $E((r_0,\theta_0), \rho)$ 
is fully contained in ${\cal D},$ the maximal admissible $\rho>0$ is given by 
$$\rho = \min_{(r,\theta) \in \partial {\cal D}}\sqrt{\left(\log \frac{r}{r_0}\right)^2 + (\theta -\theta_0)^2},$$
 where $\partial {\cal D}$ denotes  the boundary of ${\cal D}.$
\vskip0,4cm
\begin{figure}[htbp]
\centering
\includegraphics[width=50mm]{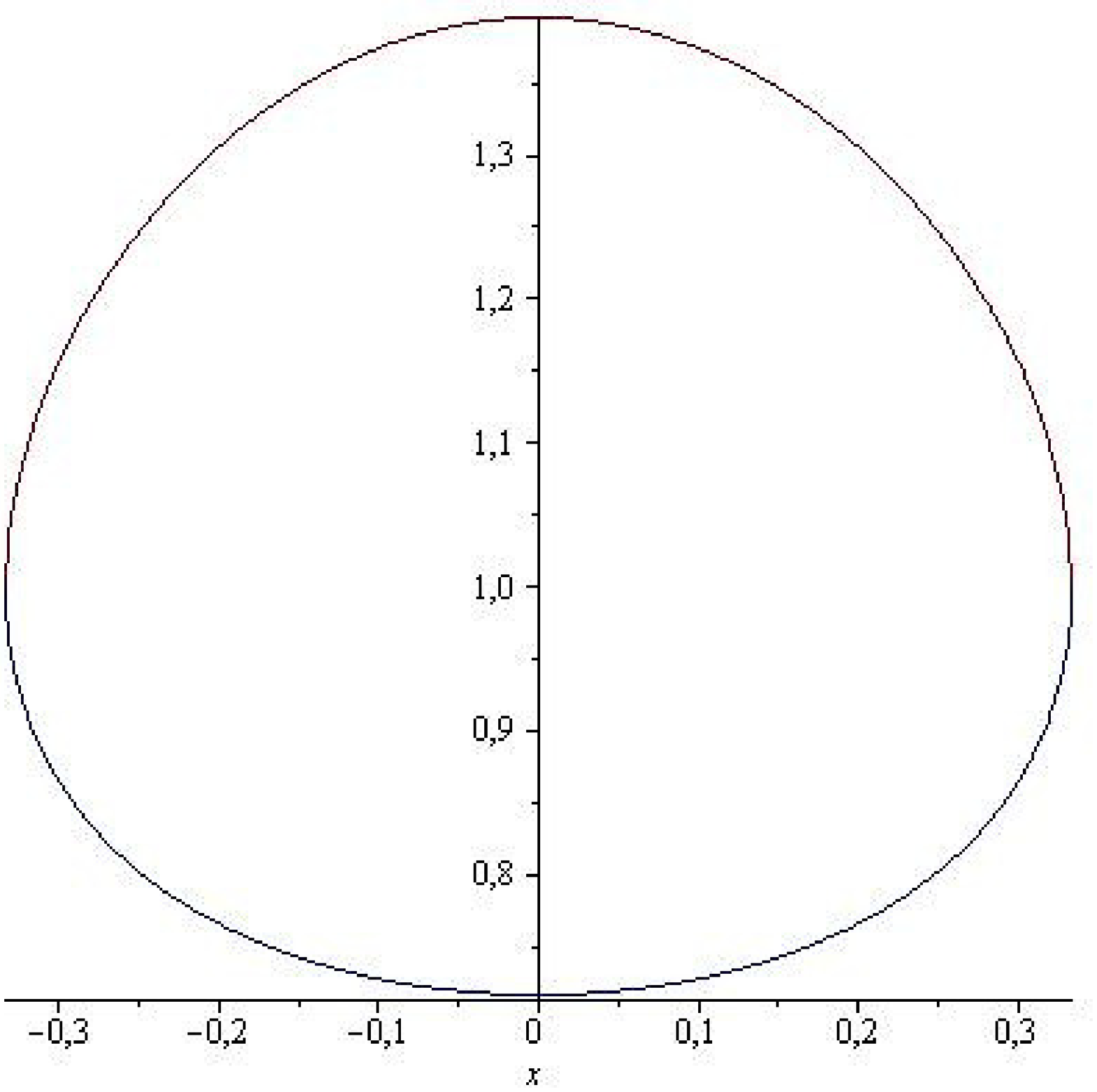}%
\qquad\qquad
\includegraphics[width=50mm]{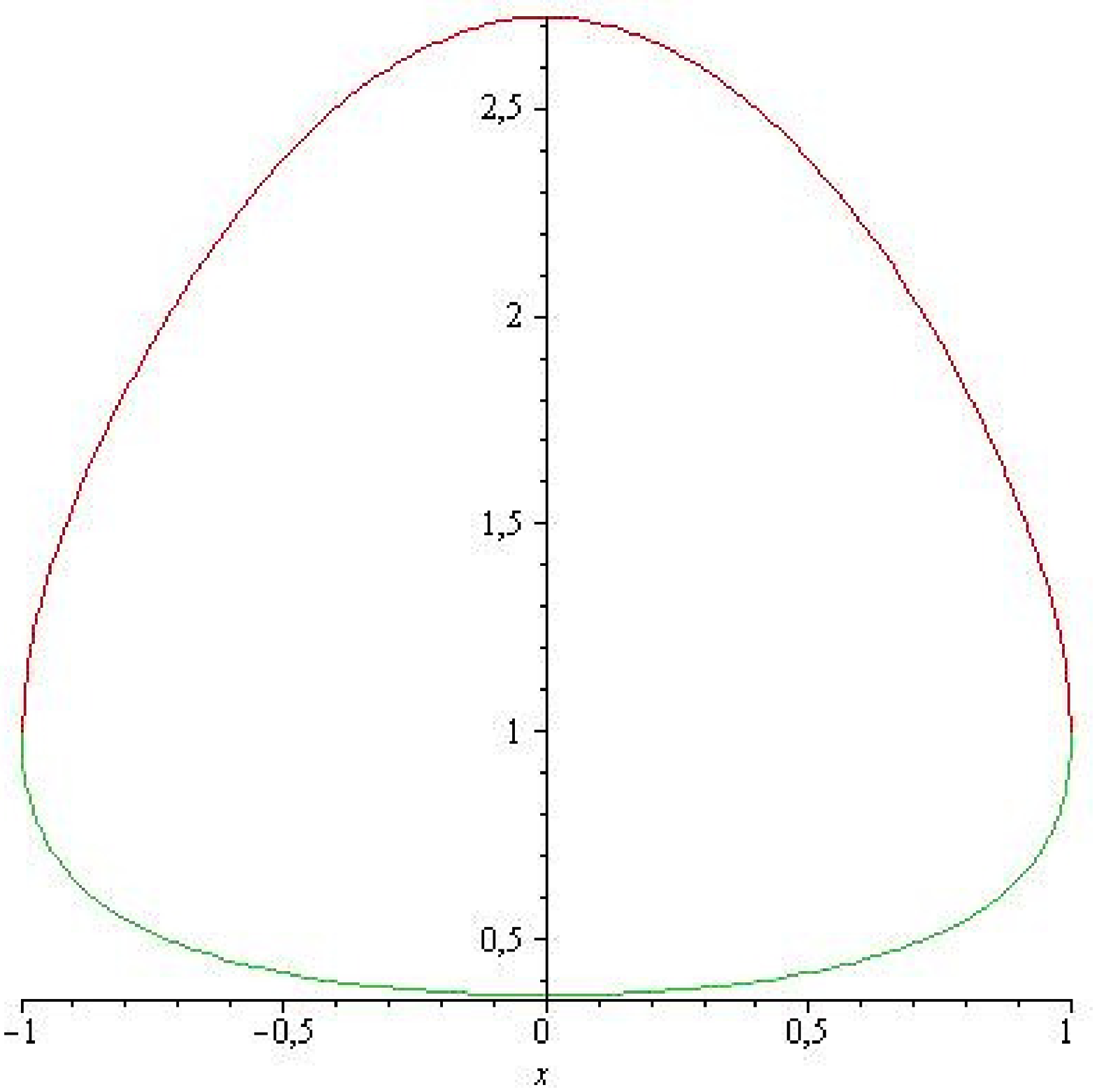}
\caption{\small The polar-disks around the point $(\theta_0, r_0)=(0,1)$ with $\rho=1/3$ and $\rho=1.$} 
\end{figure}

\vskip0,4cm
\begin{figure}[htbp]
\centering
\includegraphics[width=50mm]{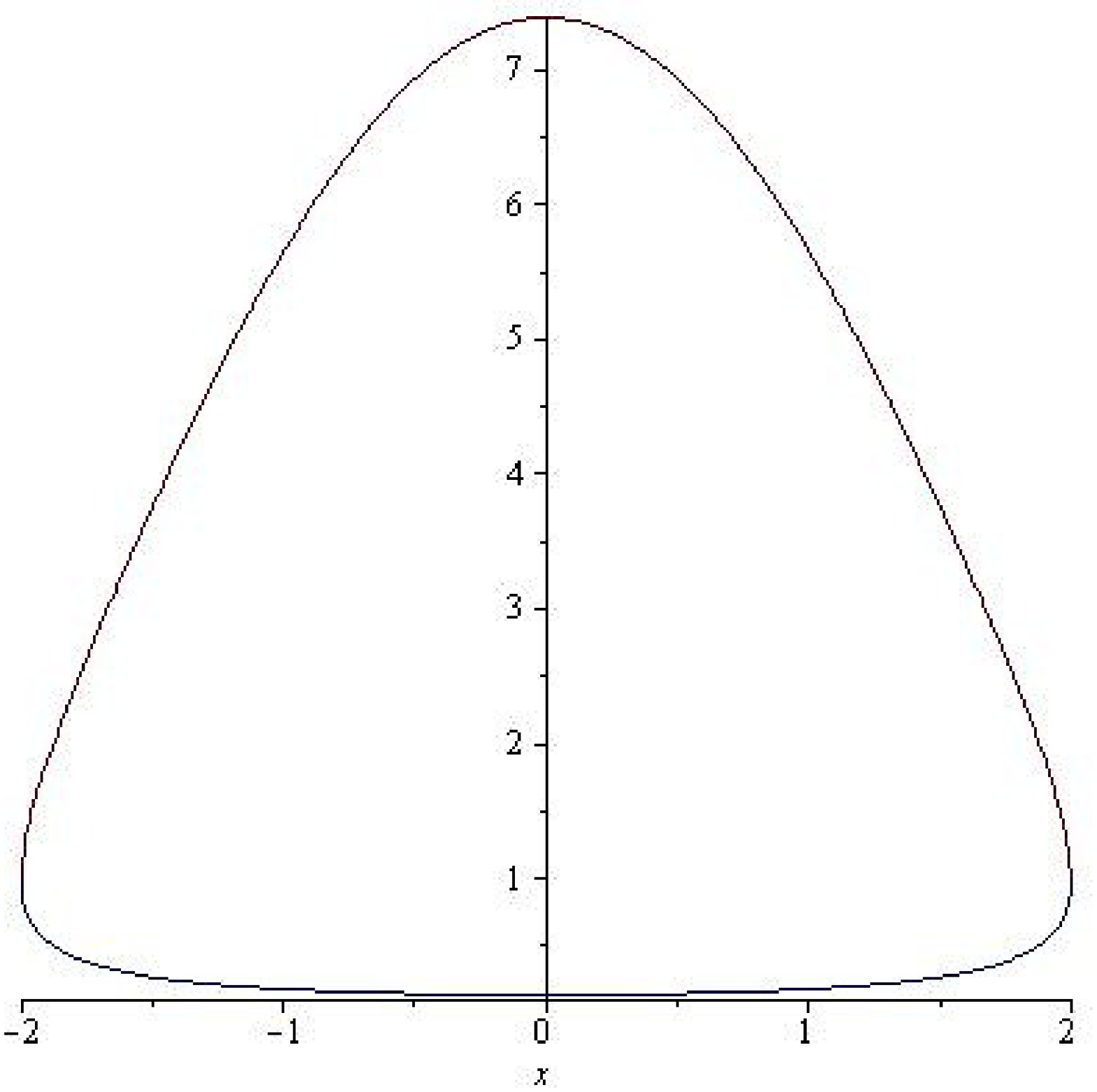}%
\qquad\qquad
\includegraphics[width=50mm]{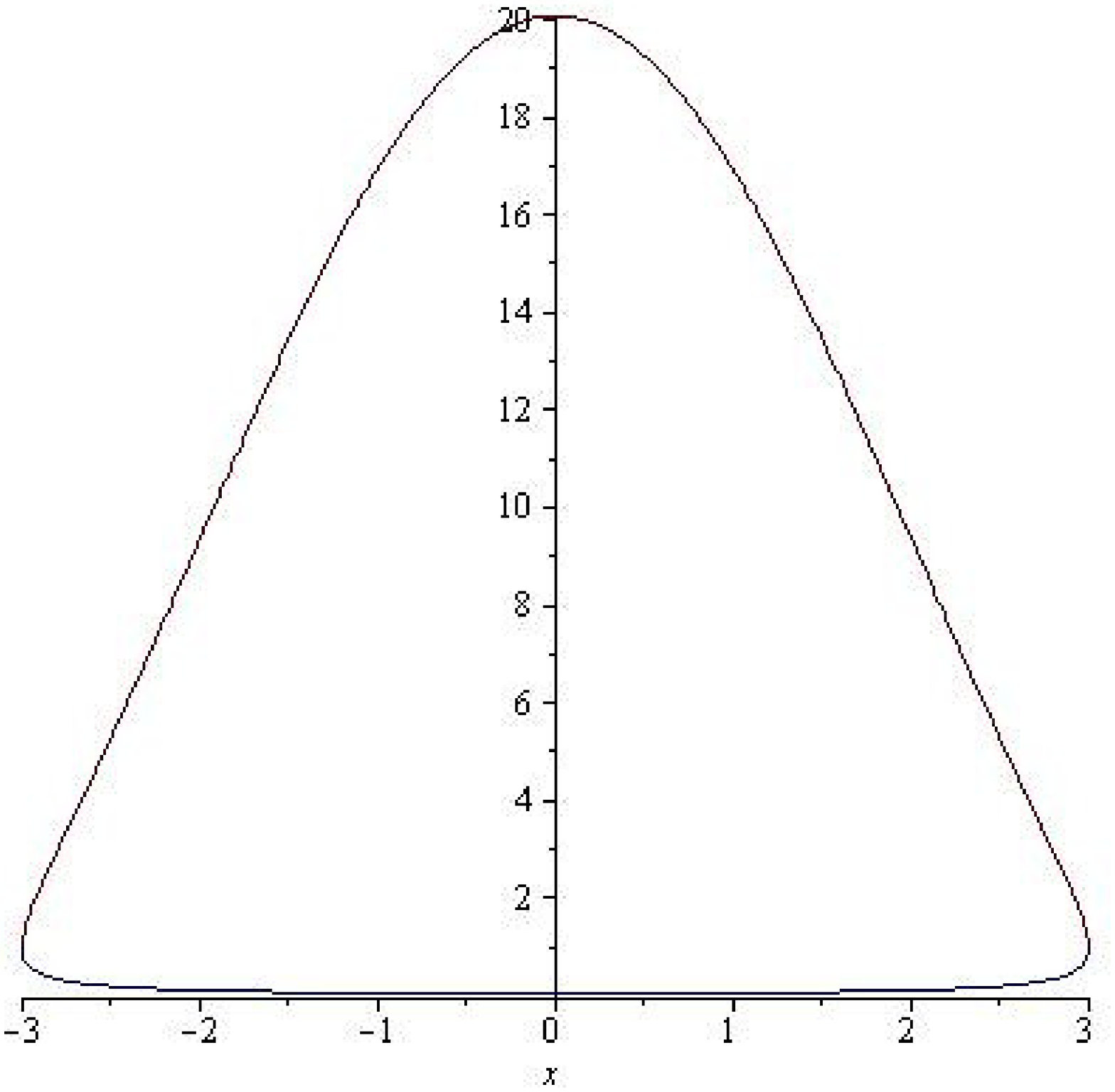}
\caption{\small The polar-disks around the point $(\theta_0, r_0)=(0,1)$ with $\rho=2$ and $\rho=3.$} 
\end{figure}

\newtheorem{Proposition}{Proposition}
\begin{Proposition}\label{taylor}
Let $f:{\cal D}\rightarrow \mathbb{C}$ be polar-analytic on ${\cal D}.$ If $(r_0, \theta_0) \in {\cal D},$ then there exists an expansion
$$f(r,\theta) = \sum_{k=0}^\infty a_k \bigg(\log\frac{r}{r_0} + i(\theta - \theta_0)\bigg)^k,$$
converging uniformly on every polar-disk $E((r_0, \theta_0), \rho) \subset {\cal D}.$
\end{Proposition}
{\bf Proof}. Let us consider the function $g(z) = g(x+iy) := f(e^x,y),$ which is defined on the domain $A \subset \mathbb{C}.$ Writing $f = u + iv$ with real-valued functions $u$ and $v,$ 
we know that the differential equations (\ref{CRE}) hold. Therefore it is easily seen that $g$ satisfies the Cauchy--Riemann equations on $A$, and so $g$ is 
an analytic function on $A.$ Hence, for $z_0 = x_0 + iy_0 \in A$ we have 
$$g(z) = \sum_{k=0}^\infty \frac{g^{(k)}(z_0)}{k!}(z-z_0)^k,$$
convergent on every disk centered in $z_0$ and contained in $A.$ Writing $z= x+iy, z_0 = x_0 + iy_0,$ and setting $r= e^x, r_0 = e^{x_0}, \theta = y, \theta_0 = y_0,$ we have, for $a_k = g^{(k)}(z_0)/k!,$
$$f(r, \theta) = \sum_{k=0}^\infty a_k \bigg(\log\frac{r}{r_0} + i(\theta -\theta_0)\bigg)^k,$$
convergent on every polar-disk $E((r_0, \theta_0), \rho)$ with sufficiently small $\rho >0.$  $\Box$
\vskip0,4cm
In particular, if $f:\mathbb{H} \rightarrow \mathbb{C}$ is polar-analytic, then the function $g$ in the proof of Proposition \ref{taylor} is an entire function. In this situation, we may choose $r_0=1, \theta_0=0$ and obtain
$$f(r,\theta) = \sum_{k=0}^\infty a_k (\log r + i\theta)^k,$$
convergent everywhere on $\mathbb{H}.$ 
\vskip0,4cm
For example, given the function $\sin z,$  which is $2\pi$-periodic with respect to $\theta,$ the function $f(r,\theta) := \sin(re^{i\theta})$ has the expansion of Proposition \ref{taylor} with $g(z) :=\sin (\exp(z))$ but also the simpler expansion
$$f(r,\theta)= \sum_{k=0}^\infty \frac{(-1)^k}{(2k+1)!}r^{2k+1}e^{i(2k+1)\theta} \quad \quad ((r,\theta) \in \mathbb{H}).$$
\vskip0,4cm

\vskip0,3cm
The following proposition on line integrals for polar-analytic functions includes an analogue of Cauchy's fundamental theorem of complex function theory and represents an important result of this new open function theory. For sake of simplicity we consider functions $f$ defined on $\mathbb{H}.$
Here a piecewise continuously differentiable curve will be called a {\it regular curve}. 
\begin{Proposition}\label{lineintegral}
Let $f$ be a polar-analytic function on $\mathbb{H}$ and let $(r_1, \theta_1)$ and $(r_2, \theta_2)$ be any two points in $\mathbb{H}.$ Then the line integral
\begin{eqnarray}\label{int}
\int_\gamma f(r,\theta) e^{i\theta}(dr + ird\theta)
\end{eqnarray}
has the same value for each regular curve $\gamma$ in $\mathbb{H}$ that starts at $(r_1,\theta_1)$ and ends $(r_2,\theta_2).$ In particular, the integral vanishes for closed regular curves.
\end{Proposition}
{\bf Proof}. Recalling equations (\ref{CRE2}), we easily verify that
$$\frac{\partial}{\partial \theta}\bigg[f(r, \theta)e^{i\theta}\bigg] = 
\frac{\partial}{\partial r}\bigg[f(r, \theta) ir e^{i\theta}\bigg].$$
By a theorem of Schwartz, this implies that the integrand in (\ref{int}) is an exact differential on $\mathbb{H}$, that is, there exists a function $F :\mathbb{H} \rightarrow \mathbb{C}$ such that
$$ \frac{\partial F}{\partial r}(r,\theta) = f(r, \theta)e^{i\theta}\quad \mbox{and} \quad 
\frac{\partial F}{\partial \theta}(r, \theta) = f(r,\theta) ir e^{i\theta},$$
and so the integral in (\ref{int}) is equal to $F(r_2,\theta_2) - F(r_1, \theta_1).$  $\Box$
\vskip0,4cm
We can also establish a converse of Proposition \ref{lineintegral} in the spirit of Morera's theorem. We recall that there exist several equivalent versions of Morera's theorem, using general closed regular curves, triangles or rectangles (see e.g. \cite[p. 208]{RU}, \cite[p. 169]{BC}, \cite[Sec. 7.2, p. 79]{BN}). We prefer a version with rectangles which seems more  suitable for our purposes.
\begin{Proposition}\label{morera}
Let $f:\mathbb{H}\rightarrow \mathbb{C}$ be a continuous function. If the integral (\ref{int}) is zero on the boundary of every rectangle $R\subset \mathbb{H},$ then $f$ is polar-analytic on $\mathbb{H}.$
\end{Proposition}
{\bf Proof}. Setting $g(z) = g(x+iy):= e^{x+iy}f(e^x,y),$ we obtain a continuous function on $\mathbb{C}.$ Moreover, if $R = [a,b]\times[c,d],$ then under the transformation $r= e^x$ and $\theta = y,$ we obtain another rectangle $R_1:=[\log a, \log b] \times [c,d]$ on the complex plane, and conversely, each rectangle on $\mathbb{C}$ is associated with a rectangle on $\mathbb{H}$ through the inverse transformation. Since
$$ \int_{\partial R}  f(r,\theta) e^{i\theta}(dr + ird\theta) = \int_{\partial R_1} g(x+iy)d(x+iy) = 0,$$
we deduce from Morera's  theorem ( \cite[Sec. 7.2, p. 79]{BN}) that the function $g$ is an entire function. Hence the function $(x,y) \longmapsto f(e^x,y)$ is also entire. Thus it satisfies the Cauchy-Riemann equations on $\mathbb{C}.$ This implies that $f$ satisfies equations (\ref{CRE}) and so the assertion follows. $\Box$

%------------------------------------------------
\section{Geometric properties of polar-analytic functions}

A remarkable geometric property of an analytic function $g$ is that at any
point $z_0$ with $g'(z_0)\ne0$ it preserves angles and orientation (see \cite[page 145]{BN}). More
precisely, this means the following. Let $\gamma_1$ and $\gamma_2$ be
smooth arcs that intersect in $z_0$ where they have tangents
$\boldsymbol{t}_1$ and $\boldsymbol{t}_2$, respectively.  Suppose that we
have to rotate $\boldsymbol{t}_1$ in the mathematical positive sense by an
angle $\alpha$ around $z_0$ in order that it coincides with
$\boldsymbol{t}_2$. Then the same is true for the arcs $g\circ\gamma_1$
and $g\circ\gamma_2$, that is, if $\boldsymbol{\tau}_1$ and
$\boldsymbol{\tau}_2$, respectively, are their tangents in $g(z_0)$, then
we have to rotate $\boldsymbol{\tau}_1$ in the mathematical positive sense
by the angle $\alpha$ around $g(z_0)$ in order that it coincides with
$\boldsymbol{\tau}_2$.
We cannot expect that this property extends to polar-analytic functions.
In the sequel we want to study  how angles may change for this class of
functions.

Let $f\,:\, (r,\theta) \mapsto u(r,\theta)+ iv(r,\theta)$ be
polar-analytic. 
Let
$$ {r\choose\theta}\,=\,{r_0\choose\theta_0} + t {c_j\choose s_j}
\qquad (t\in\mathbb{R}; \,j=1,2),$$
where
$$c_j= \cos \phi_j,\quad s_j=\sin \phi_j \qquad 
\left(\phi_j\in \left]-\frac{\pi}{2}, \frac{\pi}{2}\right];\, j=1,2\right),$$
be the tangents of two smooth arcs $\gamma_1$ and $\gamma_2$
intersecting at $(r_0,\theta_0)$, and suppose that 
$(D_{\text{pol}}f)(r_0,\theta_0)\ne0.$
By a standard formula for the angle $\alpha$ between two vectors, we have
\begin{equation}\label{4}
\cos \alpha\, =\, 
\frac{{\displaystyle \left\langle {c_1\choose s_1},{c_2\choose
s_2}\right\rangle}}%
{{\displaystyle\left\|{c_1\choose s_1}\right\|_2 \cdot
	\left\|{c_2\choose s_2}\right\|_2}},
\end{equation}
where $\langle \cdot , \cdot\rangle$ is the standard inner product in $\mathbb{R}^2$, 
and this
gives $\cos \alpha = \cos(\phi_2-\phi_1)$ in accordance with our
specification of the tangents.
Now, using (\ref{3}), we conclude that the tangents of the images of
$\gamma_1$ and $\gamma_2$ under the mapping $f$ intersect in the point
$f(r_0,\theta_0)$ under an angle $\beta$ satisfying
\begin{equation}\label{5}
\cos \beta\, =\, 
\frac{{\displaystyle \left\langle J(r_0,\theta_0){c_1\choose s_1},
	J(r_0,\theta_0){c_2\choose s_2}\right\rangle}}%
{{\displaystyle\left\|J(r_0,\theta_0){c_1\choose s_1}\right\|_2 \cdot
	\left\|J(r_0,\theta_0){c_2\choose s_2}\right\|_2}}.
\end{equation}
Considering the numerator and the denominator of the right-hand side
separately, we find by a short calculation that
\begin{align*}
 \left\langle J(r_0,\theta_0){c_1\choose s_1}, 
	J(r_0,\theta_0){c_2\choose s_2}\right\rangle &=
\left(c_1 \frac{\partial u}{\partial r} - s_1 r_0 \frac{\partial v}{\partial r}\right)
\left(c_2 \frac{\partial u}{\partial r} - s_2 r_0 \frac{\partial v}{\partial
r}\right)\\
& \quad +
\left(c_1 \frac{\partial v}{\partial r} + s_1 r_0 \frac{\partial u}{\partial r}\right)
\left(c_2 \frac{\partial v}{\partial r} - s_2 r_0 \frac{\partial u}{\partial
r}\right)\\
&= \left[\left(\frac{\partial u}{\partial r}\right)^2+
\left(\frac{\partial v}{\partial r}\right)^2\right]\left(c_1c_2 + r_0^2 s_1 s_2\right)
\end{align*}
and
\begin{align*}
\left\|J(r_0,\theta_0){c_j\choose s_j}\right\|_2^2 &=
\left(c_j \frac{\partial u}{\partial r} - s_jr_0 \frac{\partial v}{\partial
r}\right)^2+
\left(c_j \frac{\partial v}{\partial r} + s_jr_0 \frac{\partial u}{\partial
r}\right)^2\\
&= \left[\left(\frac{\partial u}{\partial r}\right)^2+
\left(\frac{\partial v}{\partial r}\right)^2\right]\left(c_j^2 + r_0^2 s_j^2\right)
\end{align*}
for $j=1, 2.$
Substituting these expressions in (\ref{5}), we see that the terms depending
on $f=u+iv$ factor out completely, and we obtain
\begin{equation}\label{6}
\cos \beta \,=\, \frac{c_1c_2+ r_0^2s_1s_2}{\sqrt{c_1^2+r_0^2s_1^2}\cdot
	\sqrt{c_2^2+r_0^2s_2^2}}
=
\frac{{\displaystyle \left\langle {c_1\choose r_0s_1},{c_2\choose
r_0s_2}\right\rangle}}%
{{\displaystyle\left\|{c_1\choose r_0s_1}\right\|_2 \cdot
	\left\|{c_2\choose r_0s_2}\right\|_2}}.
\end{equation}
This formula shows that $\cos \beta$ does not depend on $f$.
Indeed, $\cos \beta$
only depends on the local geometric situation in $\mathbb{H}$ expressed by $r_0$,
$\phi_1$ and $\phi_2$. It also does not depend on $\theta_0$. Moreover,
comparison of (\ref{4}) and (\ref{6}) shows that angles are preserved when
$r_0=1$.
These remarkable properties may be summarized as follows.
\begin{Proposition}\label{prop1}
Let $f$ be polar-analytic in a neighborhood of a point
$(r_0,\theta_0)\in\mathbb{H}$, and suppose that $(D_{\text{\rm
pol}}f)(r_0,\theta_0)\ne0.$ Let $\gamma_1$ and $\gamma_2$ be two smooth
arcs that intersect in $(r_0,\theta_0)$ under an angle $\alpha$. Then the
arcs $f\circ\gamma_1$ and $f\circ\gamma_2$ intersect in $f(r_0,\theta_0)$
under an angle $\beta$ that depends on the local geometric situation in
$\mathbb{H}$ but does not depend on $f$ and $\theta_0$. Furthermore, $\beta=\alpha$ when $r_0=1$.
\end{Proposition}

When $\phi_1, \phi_2\in ]-\pi/2, \pi/2[$, then $c_1c_2\ne0$ and therefore
we may rewrite (\ref{6}) as
\begin{equation}\label{7}
\cos \beta \,=\, \frac{1+ r_0^2t_1t_2}{\sqrt{1+r_0^2t_1^2}\cdot
	\sqrt{1+r_0^2t_2^2}},
\end{equation}
where $t_j=\tan \phi_j$ for $j=1, 2.$ For studying the dependence on
$r_0$, we denote the right-hand side of (\ref{7}) by $C(r_0^2)$ and find by
standard calculus
$$
C'(r_0^2)\left\{\begin{array}{ccc} 
		<0 & \hbox{ if } & r_0^2t_1t_2<1,\\[1ex]
		=0 & \hbox{ if } & r_0^2t_1t_2=1,\\[1ex]
		>0 & \hbox{ if } & r_0^2t_1t_2>1.
		\end{array}\right.
$$
If $t_1t_2<0$, then $\cos \beta$ is decreasing with growing $r_0$ and so
$\beta$ is increasing. We have $\beta\to0$ as $r_0\to 0$ and $\beta\to
\pi$ as $r_0\to\infty$. If $t_1t_2>0$, then $\cos \beta$ has a minimum
at $r_0=1/\sqrt{t_1t_2}$ and is decreasing to the left of that point and
increasing to the right.

Next let us consider the situation where $\phi_1=0$ and $\phi_2=\alpha\in
]0, \pi/2]$, which has not been covered so far. Then formula (\ref{6})
yields
$$ \cos \beta\,=\, \frac{c_2}{\sqrt{c_2^2+ r_0^2s_2^2}}\,=\,
\frac{\cos \alpha}{\sqrt{\cos^2\alpha+ r_0^2 \sin^2\alpha}}\,.$$
For $\alpha\in ]0, \pi/2[$, we find that $\beta$ is increasing with
growing $r_0$. We have
$$\begin{array}{ccc}
\beta \to 0 & \hbox{ as } & r_0\to 0, \\[1ex]
\beta = \alpha  & \hbox{ if } & r_0=1, \\[1ex]
\beta \to \pi/2 & \hbox{ as } & r_0\to \infty.
\end{array}
$$
For $\alpha=\pi/2$, we see that $\cos \beta =0$ for all $r_0>0$ and so
always $\beta=\pi/2.$ This is another remarkable observation which may be
stated as follows.
\begin{Proposition} \label{prop2}
Let $\mathcal{D}$ be a domain in $\mathbb{H}$, and let $f\,:\, \mathcal{D} \to\mathbb{C}$
be polar-analytic with $D_{\text{\rm pol}}f \ne0$. Denote by $\mathcal{N}$ an
orthogonal net of lines parallel to the axes of $\mathbb{H}$. Then $f$ maps
$\mathcal{N}\cap \mathcal{D}$ onto an orthogonal net of curves in
$f(\mathcal{D})$.
\end{Proposition}

As an example, consider the polar-analytic function $f$ given by
$f(r,\theta):=e^{-i\theta}/r.$ Let $(a_n)_{n\in\mathbb{N}}$ be a sequence of
positive numbers and let $(b_n)_{n\in\mathbb{N}}$ be a sequence of real numbers.
Then the lines
$$ L_n\,:=\, \left\{(r,\theta)\in\mathbb{H}\::\: r=a_n\right\}$$
are mapped by $f$ to circles of radius $1/a_n$ centered at $0$. The half-lines
$$M_n\,:=\, \left\{(r,\theta)\in\mathbb{H}\::\: \theta=b_n\right\}$$
are mapped by $f$ to the rays $\{te^{-ib_n}\,:\,t>0\}$.
Clearly, the circles and the rays constitute an orthogonal net in $\mathbb{C}$.
However, this example reveals that, other than in the case of an analytic
function, orientation is {\it not} preserved.

%----------------------------------------------------------------
\section{Why we use polar-analytic functions?}

So far, the chief applications of the concept of polar analyticity have been established in Mellin analysis and in the realm of quadrature formulae. For the first, one can obtain an extension of the classical Paley-Wiener theorem of Fourier analysis to the setting of Mellin transforms. Indeed, the notion of a polar-analytic function gives a simple definition of the so-called Mellin-Bernstein spaces, without resorting to Riemann surfaces and analytical branches. These spaces characterize, in a precise sense, the functions $f: \mathbb{R}^+ \rightarrow \mathbb{C}$ with compactly supported Mellin transform (see \cite{BBMS1}, \cite{BBMS2}). Moreover, polar analyticity enables us to define in a simple way the Hardy spaces in Mellin frame (see \cite{BBMS1}). 
In Mellin frame, the Mellin derivative, for functions $f: \mathbb{R}^+ \rightarrow \mathbb{C},$ introduced in \cite{BJ}, is defined, for $c \in \mathbb{R},$ by 
$$\Theta_cf(x) := xf'(x) + cf(x) \quad (x>0),$$
provided that $f'(x)$ exists.

Note that in case $c=0$ we obtain the original definition of the Mellin derivative as given in \cite{MA}, also known as "Euler" derivative due to its use in the theory of Euler type (partial) differential equations, which are linear equations with variable (polynomial) coefficients (see e.g. the recent paper \cite{DL}).

 For functions defined in the polar plane $\mathbb{H},$ we have the corresponding notion of Mellin polar-analytic function and the Mellin polar-derivative on putting (see \cite{BBMS1})
$$\widetilde{\Theta}_cf(r,\theta):= re^{i\theta}(D_{{\rm pol}} f)(r, \theta) + cf(r, \theta).$$
We can see that if $\varphi(\cdot):= f(\cdot, 0),$ then $\widetilde{\Theta}_c \varphi(r) = \Theta_c\varphi(r),$ for $r>0.$

Concerning quadrature formulae, polar analyticity gives a characterization of certain function spaces of Hardy type in terms of the speed of the convergence to zero of the remainders in certain Gaussian quadrature formulae for functions defined on the positive real axis (see \cite{BBMS3}).


\begin{thebibliography}{99}
\small
\bibitem{AH} L. V. Ahlfors, Complex Analysis, McGraw-Hill Int. Eds, Third Edition, 1979.
\bibitem{BN} J. Bak and D.J. Newman, Complex Analysis, Springer-Verlag, New York, 1982.
\bibitem{BBMS1} C. Bardaro, P.L. Butzer, I. Mantellini and G. Schmeisser, A fresh approach to the Paley-Wiener theorem for Mellin transforms and the 
	Mellin-Hardy spaces,  Math. Nachr. 290, (2017), 2759--2774.
\bibitem{BBMS2} C. Bardaro, P.L. Butzer, I. Mantellini and G. Schmeisser, A generalization of the Paley--Wiener theorem for Mellin transforms and 
	metric characterization of function spaces, Frac. Calc. Appl. Anal., 20(5), (2017), 1216--1238.
\bibitem{BBMS3} C. Bardaro, P.L. Butzer, I. Mantellini and G. Schmeisser, Quadrature formulae for the positive real axis in the setting of Mellin analysis: sharp error estimates in terms of the Mellin distance, submitted for publication, December 2017. arXiv: 1802.03952v1.
\bibitem{BBMS4} C. Bardaro, P.L. Butzer, I. Mantellini and G. Schmeisser, Development of a complex function theory upon a  new concept of polar-analytic functions, March 2018. arXiv: 1803.04258v1.
\bibitem{BC} J.W. Brown and R.V. Churchill, Complex Variables and Applications, Eighth Ed., McGraw-Hill, 2009.
\bibitem{BJ} P.L. Butzer and S. Jansche, A direct approach to the Mellin transform, J. Fourier Anal. Appl., 3, (1997), 325--375.
\bibitem{DL} P. Doma\'nski and M. Langenbruch, Interpolation of holomorphic functions and surjectivity of Taylor coefficients multipliers, Adv. Math., 293 (2016), 782--855
\bibitem{HI} E.V. Hille, Analytic Function Theory, vol I, Chelsea Publ. Co., N.Y., 1959.
\bibitem{KA1} T. Kaiser, The Dirichlet problem in the plane with semianalytic raw data, quasi analyticity, and $o$-minimal structure, Duke Math. J., 147(2), (2009), 285--314
\bibitem{KA2} T. Kaiser, The Riemann mapping theorem for semianalytic domains and $o$-minimality, Proc. London Math. Soc. 98(3), (2009), 427--444.
\bibitem{KRS} T. Kaiser, J.-P. Rolin and P. Speissegger, Transition maps at non-resonat hyperbolic singularities are $o$-minimal, J. reine angew. Math. 636, (2009), 1--45.
\bibitem{MA} R.G. Mamedov, The Mellin transform and approximation theory, Elm (Baku), 1991.
\bibitem{RU} W. Rudin, Real and Complex Analysis, McGraw-Hill, Third Edition, New York, 1987.
\bibitem{SI} H. Silverman, Polar form of Cauchy-Riemann equations, PRIMUS, 10(3), (2000), 241--245.


\end{thebibliography}
\end{document}